\documentclass[12pt,reqno]{amsart}

\usepackage{amsthm}
\usepackage{enumerate}
\usepackage[numbers,sort&compress]{natbib}
\usepackage{url}
\usepackage{pstricks,pst-plot,pst-node}

\newcommand{\reals}{\mathbb R}          
\newcommand{\nats}{\mathbb N}           
\newcommand{\dist}{\mathrm{dist}}

\newcommand{\SLO}{SLO}
\newcommand{\SLOst}{SLO$^\ast$\/}

\theoremstyle{plain}
\newtheorem{thm}{Theorem}
\newtheorem{lemma}{Lemma}
\newtheorem{proposition}[lemma]{Proposition}
\newtheorem{cor}[thm]{Corollary}

\theoremstyle{definition}
\newtheorem{definition}{Definition}

\theoremstyle{remark}
\newtheorem*{problem}{Problem}
\newtheorem*{remark}{Remark}

\newcommand{\vertex}[2]{\cnode*(#2){4pt}{#1}}
\newcommand{\bdvertex}[2]{\cnode[fillstyle=solid,fillcolor=white](#2){4pt}{#1}}
\newcommand{\edge}[2]{\ncline[nodesep=0pt]{-}{#1}{#2}}
\newcommand{\edgedash}[2]{\ncline[nodesep=0pt,linestyle=dashed,dash=4pt 2pt]{-}{#1}{#2}}

\usepackage{algorithmic}

\begin{document}


\title{Faber-Krahn Type Inequalities for Trees}

\author{T\"urker B{\i}y{\i}ko\u{g}lu \and Josef Leydold}

\address{Max-Planck-Institute for Mathematics in the Sciences,
  Inselstra{\ss}e 22, D-04103 Leipzig, Germany}

\email{tuerker@statistik.wu-wien.ac.at}

\address{University of Economics and Business Administration,
  Department for Applied Statistics and Data Processing,
  Augasse 2-6, A-1090 Vienna, Austria}

\email{Josef.Leydold@statistik.wu-wien.ac.at}

\thanks{This work was supported by the Austrian Science Foundation (FWF),
  project no.~14094-MAT}

                  

\begin{abstract}
  The Faber-Krahn theorem states that among all boun\-ded
  domains with the same volume in $\reals^n$ (with the standard
  Euclidean metric), a ball that has lowest first Dirichlet eigenvalue.
  Recently it has been shown that a similar result holds for
  (semi-)regular trees. In this article we show that such a theorem
  also hold for other classes of (not necessarily non-regular) trees. 
  However, for these new results no couterparts in the world
  of the Laplace-Beltrami-operator on manifolds are known.
\end{abstract}

\maketitle


\markboth{T.~B{\i}y{\i}ko\u{g}lu and J.~Leydold}{%
  Faber-Krahn Type Inequalities for Trees}


\section{Introduction}

In recent years the eigenvectors of the graph Laplacian has received
increasing attention. While its eigenvalues has been investigated for
fifty years 
\citep[see e.g.][]{Biggs:1994a,Cvetkovic:1980a,Cvetkovic:1988a},
there is little known about the eigenvectors. The graph Laplacian can
be seen as the discrete analogon of the continuous
Laplace-Beltrami-operator on manifolds. When using an appropriate
definition for the gradiant on a graph rules similar to the classical
Laplace operator can be formulated, e.g.\ Green's formula.
During the last years some results for eigenfunctions of the
Laplace-Beltrami-operator have been shown to hold also for
eigenvectors of the graph Laplacian;
for example Cheeger-type inequalities \citep{deVerdiere:1993b} or nodal
domain theorems \citep{Davies;etal:2001a} exist.
However, it has turned out that there are small but subtle differences
between the discrete and the continuous case.

The Faber-Krahn inequality is another well-known result. It states
that among all bounded domains with the same volume in $\reals^n$
(with the standard Euclidean metric), a ball that
has lowest first Dirichlet eigenvalue \citep{Chavel:1984a}.
\citet{Friedman:1993a} described the idea of ``graph with boundary''
(see below). With this concept he was able to formulate Dirichlet and
Neumann eigenvalue problems for graphs. He also conjectured an
analogon to the Faber-Krahn inequality 
for regular trees. Amazingly Friedman's conjecture is false, i.e.\ in
general these trees are similar but not equal to ``balls'', see
\citep{Pruss:1998a,Leydold:1997a} for counterexamples and 
\citep{Leydold:2002a} for a statement of the result.
This example (as well as the nodal domain theorem where also some
wrong conjuctures exist, see \citep{Davies;etal:2001a}) shows that
there is much more structure in graphs than in manifolds. Conclusions
from this fact are twofold: First, some care is necessary since one's
intuition, trained on manifolds, may lead to wrong conjectures. On the
other hand we can use the opportunity to go further and try to find
these new structures. It is this second conclusion that motivates this
paper. We want to leave the world of regular graphs and look what
happens when we drop this regularity assumption.

In this article we want to formulate Faber-Krahn type theorems for
trees which need not to be regular any more. Analogous results for the
Laplace-Beltrami-operators on manifolds with non-constant curvature 
are rare.


\section{Discrete Dirichlet operator and Faber-Krahn property}

Let $G(V,E)$ be a simple (finite) undirected graph with vertex set $V$ and edge
set $E$. The \emph{Laplacian} of $G$ is the matrix 
\begin{equation}
  \Delta(G) = D(G) - A(G), 
\end{equation}
where $A(G)$ denotes the adjacency matrix of the graph and $D(G)$ is
the diagonal matrix whose entries are the vertex degrees, i.e.,
$D_{vv} = d_v$, where $d_v$ denotes the degree of vertex $v$.
We write $\Delta$ for short if there is no risk of confussion.
To state a Faber-Krahn type inequality we need a Dirichlet operator
which itself requires the notion of a boundary of a graph.

A \emph{graph with boundary} $G(V_0\cup\partial V,E_0\cup\partial E)$
consists of a set of interior vertices $V_0$, boundary vertices
$\partial V$, interior edges $E_0$ that connect interior vertices,
and boundary edges $\partial E$ that join interior vertices with
boundary vertices \citep{Friedman:1993a}.
There are no edges between two boundary vertices.

In the following we assume that every boundary vertex has degree $1$
and every interior vertex has degree at least 2, i.e.\
a vertex is a boundary vertex if and only if it has degree 1.
We also assume that both the set interior vertices $V_0$ and the
set of boundary vertices $\partial V$ are not empty.
Balls are of particular interest for our investigations.
A \emph{ball} $B(v_0,r)$ with center $v_0$ and radius $r\in\nats$ is
a connected graph where every boundary vertex $w$ has geodesic
distance $\dist(v_0,w)=r$.

A discrete Dirichlet operator is the graph Laplacian $\Delta$
which acts on vectors only that vanish in all boundary vertices. For a
motivation of this definition see \citep{Friedman:1993a}.

\begin{definition}
  A \emph{discrete Dirichlet operator} $\Delta_0$ is the graph
  Laplacian restricted to interior vertices, i.e.
  \begin{equation}
    \Delta_0 = D_0 - A_0
  \end{equation}
  where $A_0$ is the adjacency matrix of the graph induced by the
  interior vertices, $G(V_0,E_0)$, and
  where $D_0$ is the degree matrix with the vertex degrees 
  in the whole graph $G(V_0\cup\partial V,E_0\cup\partial E)$ as its
  entries. 
\end{definition}

Notice that $\Delta_0$ is obtained from the graph Laplacian $\Delta$
by deleting all rows and columns that correspond to boundary vertices.
Thus any edges between two boundary vertices have no influence on the
Dirichlet operator. Thus we have eliminated such edges by definition
for the sake of simplicity.

\begin{definition}[Faber-Krahn property]
  We say that a graph with boundary has the \emph{Faber-Krahn
    property} if it has lowest first Dirichlet eigenvalue among all
  graphs with the same ``volume'' in a particular graph class.
\end{definition}

This definition raises two questions:
(1) What is the ``volume'' of a graph, and (2) what is an appropriate
graph class (besides the trivial requirement that it must contain the
graph $G$ in question)?

\begin{figure}[ht]
  \centering
  {
\newcommand{\grayedge}[2]{\ncline[linecolor=gray,nodesep=0pt,linestyle=solid]{-}{#1}{#2}}
\begin{pspicture}(-3.3,-3.3)(3.3,3.3)
  \psccurve[linecolor=gray,linewidth=1.5pt,linestyle=dashed,dash=4pt 2pt]
  (1.73205,1.)(0.258819,2.29926)(-0.38823,2.448899)(-1.99156,1.74755)
  (-2.31492,0.88823)(-1.4240,-0.9665)(-0.25,-1.4330)(1.3662,-0.5)
  \vertex{0}{0,0}
  \vertex{1}{1,0}
  \vertex{2}{-0.5,0.866025}
  \vertex{3}{-0.5,-0.866025}
  \edge{0}{1}
  \edge{0}{2}
  \edge{0}{3}
  \bdvertex{4}{1.73205,1.}
  \vertex{5}{0,2.}
  \vertex{6}{-1.73205,1.}
  \pnode(-1.73205,-1.){7}
  \pnode(0.,-2.){8}
  \pnode(1.73205,-1.){9}
  \edge{1}{4}
  \grayedge{1}{9}
  \edge{2}{5}
  \edge{2}{6}
  \grayedge{3}{7}
  \grayedge{3}{8}
  \pnode(2.89778,0.776457){10}
  \pnode(2.12132,2.12132){11}
  \pnode(0.776457,2.89778){12}
  \pnode(-0.776457,2.89778){13}
  \pnode(-2.12132,2.12132){14}
  \pnode(-2.89778,0.776457){15}
  \pnode(-2.89778,-0.776457){16}
  \pnode(-2.12132,-2.12132){17}
  \pnode(-0.776457,-2.89778){18}
  \pnode(0.776457,-2.89778){19}
  \pnode(2.12132,-2.12132){20}
  \pnode(2.89778,-0.776457){21}
  \grayedge{4}{10}
  \grayedge{4}{11}
  \grayedge{5}{12}
  \grayedge{5}{13}
  \grayedge{6}{14}
  \grayedge{6}{15}
  \grayedge{7}{16}
  \grayedge{7}{17}
  \grayedge{8}{18}
  \grayedge{8}{19}
  \grayedge{9}{20}
  \grayedge{9}{21}
  \bdvertex{b37}{-1.4240,-0.9665}
  \edge{3}{b37}
  \bdvertex{b38}{-0.25,-1.4330}
  \edge{3}{b38}
  \bdvertex{b19}{1.3662,-0.5}
  \edge{1}{b19}
  \bdvertex{b512}{0.258819,2.29926}
  \edge{5}{b512}
  \bdvertex{b513}{-0.38823,2.448899}
  \edge{5}{b513}
  \bdvertex{b614}{-1.99156,1.74755}
  \edge{6}{b614}
  \bdvertex{b615}{-2.31492,0.88823}
  \edge{6}{b615}
\end{pspicture}
}
  \caption{The class of trees considered by \citet{Friedman:1993a}
    can be obtained by cutting connected subsets out of the geometric
    representation of an infinite $d$-regular tree.\hfill\break
    ($\bullet$ \ldots\ interior vertices, $\circ$ \ldots\ boundary vertices)}
  \label{fig:infinite-tree}
\end{figure}
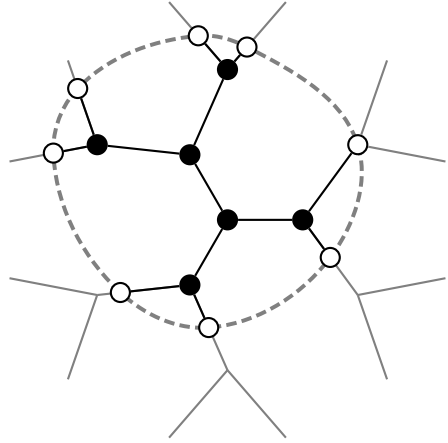

\citet{Pruss:1998a} used the number of edges of an unweighted tree
as volume and the class of semi-$d$-regular trees with boundary. In
such a tree every interior vertex has the same degree $d$ whereas
every boundary vertex has degree 1. This idea can be extended to
weighted trees \citep{Friedman:1993a}, where edge weights are
represented by the reciprocal lengths of arcs in a geometric
representation of the tree. The volume is then defined as the
sum of all the arc lengths of the geometric representation.
\citet{Friedman:1993a} looked at the class of all
trees, where the interior vertices have the same degree $d$, all
interior vertices have length (weight) 1 and all boundary vertices
have length at most 1. Such graphs can be obtained by cutting
out a subset of the geometric reperesentation of an infinite
(unweighted) $d$-regular tree, see Fig.~\ref{fig:infinite-tree}.

In this article we want to formulate Faber-Krahn type theorems for
(non-regular) trees. Analogous results for the
Laplace-Beltrami-operators on manifolds with non-constant curvature
are rare.
When we generalize the Faber-Krahn type theorems to arbitrary
trees, the picture of cutting out a graph fails. Instead we have to 
solve the following problem.

\begin{problem}
  Given a class $\mathcal{C}$ of graphs, where all graphs have the same
  ``volume''. Now characterize all graphs in $\mathcal{C}$ with the
  Faber-Krahn property, i.e., which minimize the first Dirichlet
  eigenvalue.
\end{problem}

Making the graph class $\mathcal{C}$ too large leads to quite simple
(non-interesting) graphs. For example, if $\mathcal{C}$ is the set of
all connected graphs with a given number of vertices as the ``volume'' 
of the graph, then graphs with the Faber-Krahn property are paths with
one terminating triangle \citep{Katsuda;Urakawa:1999a}.
If we restrict this class to trees, then we arrive at simple pathes
\citep{Katsuda;Urakawa:1999a,Katsuda;Urakawa:1998a}.
(To be precise \citet{Katsuda;Urakawa:1999a} used the more general
``non-separation property''.)

It seems natural to use the number of vertices as measure for
the ``volume'' of a graph. (Notice that this is equivalent to use the
number of edges for an unweighted tree.)
Moreover, we will consider only graph classes where both the total
number of interior vertices, $|V_0|$, and boundary vertices,
$|\partial V|$, is fixed. (For semiregular trees this is
always the case when we fix the total number of vertices.)
We will drop this requirement at the end of this article and state
some additional results in Sect.~\ref{sec:further-results}.
Hence we will look at the following classes of graphs with boundaries: 
\begin{align}
  \mathcal{T}^{(n,k)} 
  &= \{\text{$G$ is a tree, with $|V|=n$ and $|V_0|=k$}\} \\
  \mathcal{T}_d^{(n,k)} 
  &= \{G\in\mathcal{T}^{(n,k)}\colon d_v\geq d \text{ for all }v\in V_0\}
\end{align}
As it is clear that we always look at a particular class
$\mathcal{T}^{(n,k)}$ or $\mathcal{T}_d^{(n,k)}$ we will write
$\mathcal{T}$ and $\mathcal{T}_d$ for short; $n$ and $k$ have then to
be selected accordingly. We always assume that $1\leq k\leq n-1$.

Another interesting class is based on so called degree sequences.
A sequence $\pi=(d_0,\ldots,d_{n-1})$ of nonnegative integers is called 
\emph{degree sequence} if there exists a graph $G$ with $n$ vertices
for which $d_0,\ldots,d_{n-1}$ are the degrees of its vertices. 
For trees the following characterization exists.
\begin{lemma}[\citep{Harary:1969a}]
  \label{lem:treesequence}
  A degree sequence $\pi=(d_0,\ldots,d_{n-1})$ is a tree sequence 
  (i.e.\ a degree sequence of some tree) if and only if 
  every $d_i>0$ and $\sum_{i=0}^{n-1} d_i=2\,(n-1)$.
\end{lemma} 
Using this notion we can introduce another interesting graph class 
for which we want to formulate a Faber-Krahn like theorem,
\begin{equation}
  \mathcal{T}_\pi 
  = \{\text{$G$ is a tree with boundary with degree sequence $\pi$}\}\,.
\end{equation}
Notice that for a particular degree sequence $\pi$ we have
\begin{equation}
  \mathcal{T}_\pi \subseteq \mathcal{T}_{d_\pi}
  \subseteq \mathcal{T}_2 = \mathcal{T}
\end{equation}
where $d_\pi$ is the minimal degree for interior vertices of in the
degree sequence $\pi$. 

For class $\mathcal{T}$ of all trees we find a simple structure for
graphs with the Faber-Krahn property.

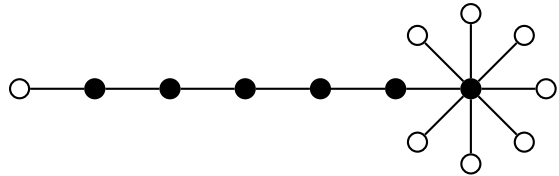
\begin{figure}[ht]
  \centering
  \begin{pspicture}(0,-1)(7,1)
  \bdvertex{A}{0,0}
  \vertex{B}{1,0}
  \vertex{C}{2,0}
  \vertex{D}{3,0}
  \vertex{E}{4,0}
  \vertex{F}{5,0}
  \vertex{G}{6,0}
  \bdvertex{I}{5.29,0.71}
  \bdvertex{J}{6,1}
  \bdvertex{K}{6.71,0.71}
  \bdvertex{H}{7,0}
  \bdvertex{L}{6.71,-0.71}
  \bdvertex{M}{6,-1}
  \bdvertex{N}{5.29,-0.71}
  \edge{A}{B}
  \edge{B}{C}
  \edge{C}{D}
  \edge{D}{E}
  \edge{E}{F}
  \edge{F}{G}
  \edge{G}{H}
  \edge{G}{I}
  \edge{G}{J}
  \edge{G}{K}
  \edge{G}{L}
  \edge{G}{M}
  \edge{G}{N}
\end{pspicture}
  \caption{A comet has the Faber-Krahn property in class
    $\mathcal{T}$. It consists of a star with diameter 2 and a path
    attached to it.\hfill\break
    ($\bullet$ \ldots\ interior vertices, $\circ$ \ldots\ boundary vertices)}
  \label{fig:klobuerstel}
\end{figure}

\begin{thm}[Klob\"ur\v{s}teltheorem]
  \label{thm:klobuerstel}
  A tree $G$ has the Faber-Krahn property in a class
  $\mathcal{T}$ if and only if G is a star with a long tail, i.e.\ a
  comet, see Fig.~\ref{fig:klobuerstel}.
  $G$ is then uniquely determined up to isomorphism.
\end{thm}

Graphs with the Faber-Krahn property in $\mathcal{T}_d$ or
$\mathcal{T}_\pi$ have a richer structure. For its description we
need additional notions.
For a tree $G$ with root $v_0$ the \emph{height}
$h(v)$ of a vertex $v$ is defined by $h(v)=\dist(v,v_0)$.
For two adjacent vertices $v$ and $w$ with $h(w)=h(v)+1$ we call 
$v$ the \emph{parent} of $w$, and $w$ a \emph{child} of $v$.
Notice that every vertex $v\not=v_0$ has exactly one parent, and
every interior vertex $w$ has at least one child vertex.

\begin{definition}[SLO-ordering]
  Let $G(V_0\cup\partial V,E_0\cup\partial E)$ be a tree
  with boundary with root $v_0$. Then a well-ordering 
  $\prec$ of the vertices is called \emph{spiral-like}
  (\emph{SLO}-ordering for short) if the following holds
  for all vertices $v,v_1,v_2,w,w_1,w_2\in V$:
  \begin{enumerate}[(S4)]
  \item[(S1)] $v\prec w$ implies $h(v)\leq h(w)$;
  \item[(S2)] if $v_1\prec v_2$ then for all children $w_1$ of $v_1$ and
    all children $w_2$ of $v_2$, $w_1\prec w_2$; 
  \item[(S3)] if $v\prec w$ and $v\in\partial V$, then $w\in\partial V$.
  \end{enumerate}
  It is called \emph{spiral-like with increasing degrees}
  (\emph{\SLOst}-ordering for short) if additonally the following
  holds 
  \begin{enumerate}[(S4)]
  \item[(S4)] if $v\prec w$ for interior vertices $v,w\in V_0$, then
    $d_v\leq d_w$.
  \end{enumerate}
  We call trees that have a SLO- or \SLOst-ordering of its
  vertices \emph{SLO-trees} and \emph{\SLOst-trees}, respectively.
\end{definition}

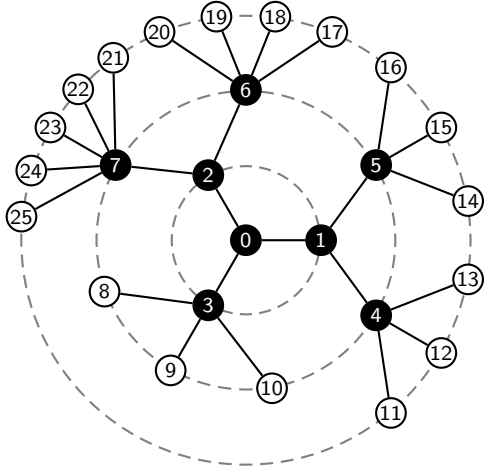
\begin{figure}[ht]
  \centering
  \newcommand{\lbbdvertex}[2]{\cnode[fillstyle=solid,fillcolor=white](#2){6pt}{#1}\rput(#2){\tiny\textsf{#1}}}
\newcommand{\lbvertex}[2]{\cnode*(#2){6pt}{#1}\rput(#2){\white\tiny\textsf{#1}}}
\begin{pspicture}(-3.3,-3.3)(3.3,3.3)
  \pscircle[linecolor=gray,linestyle=dashed](0,0){1}
  \pscircle[linecolor=gray,linestyle=dashed](0,0){2}
  \pscircle[linecolor=gray,linestyle=dashed](0,0){3}
  \lbvertex{0}{0,0}
  \lbvertex{1}{1,0}
  \lbvertex{2}{-0.5,0.866025}
  \lbvertex{3}{-0.5,-0.866025}
  \lbvertex{4}{1.73205,-1.}
  \lbvertex{5}{1.73205,1.}
  \lbvertex{6}{0,2.}
  \lbvertex{7}{-1.73205,1.}
  \lbbdvertex{8}{-1.87939,-0.68404}
  \lbbdvertex{9}{-1.,-1.73205}
  \lbbdvertex{10}{0.347296,-1.96962}
  \lbbdvertex{11}{1.92836,-2.29813}
  \lbbdvertex{12}{2.59808,-1.5}
  \lbbdvertex{13}{2.95442,-0.520945}
  \lbbdvertex{14}{2.95442,0.520945}
  \lbbdvertex{15}{2.59808,1.5}
  \lbbdvertex{16}{1.92836,2.29813}
  \lbbdvertex{17}{1.14805,2.77164}
  \lbbdvertex{18}{0.391579,2.97433}
  \lbbdvertex{19}{-0.391579,2.97433}
  \lbbdvertex{20}{-1.14805,2.77164}
  \lbbdvertex{21}{-1.76336,2.42705}
  \lbbdvertex{22}{-2.22943,2.00739}
  \lbbdvertex{23}{-2.59808,1.5}
  \lbbdvertex{24}{-2.85317,0.927051}
  \lbbdvertex{25}{-2.98357,0.313585}
  \edge{0}{1}
  \edge{0}{2}
  \edge{0}{3}
  \edge{1}{4}
  \edge{1}{5}
  \edge{2}{6}
  \edge{2}{7}
  \edge{3}{8}
  \edge{3}{9}
  \edge{3}{10}
  \edge{4}{11}
  \edge{4}{12}
  \edge{4}{13}
  \edge{5}{14}
  \edge{5}{15}
  \edge{5}{16}
  \edge{6}{17}
  \edge{6}{18}
  \edge{6}{19}
  \edge{6}{20}
  \edge{7}{21}
  \edge{7}{22}
  \edge{7}{23}
  \edge{7}{24}
  \edge{7}{25}
\end{pspicture}
  \caption{A \SLOst-tree with 8 interior and 18 boundary vertices.
    The \SLOst-ordering $\prec$ is indicated by numbers.
    Degree sequence $\pi=(3,3,3,4,4,4,5,6,1,1,\ldots,1)$.}
  \label{fig:slost-tree}
\end{figure}

Notice that \SLO-trees are almost balls, that is there exists a radius
$r$ such that $\dist(v,v_0)\in\{r,r+1\}$ for all boundary vertices
$v\in\partial V$, see Fig.~\ref{fig:slost-tree} for an example.
With this concept we can formulate Faber-Krahn type theorems for the
other graph classes, $\mathcal{T}_d$ and $\mathcal{T}_\pi$.

\begin{thm}
  \label{thm:Tdeg}
  A graph $G$ has the Faber-Krahn property in a class $\mathcal{T}_d$
  if and only if it is a \SLOst-tree 
  where exactly one interior vertex has degree 
  $d^\circ=d+\sum_{v\in V_0}(d_v-d)$
  and all other interior vertices have degree $d$.
  $G$ is then uniquely determined up to isomorphism.
\end{thm}

\begin{thm}
  \label{thm:Tdegseq}
  A graph $G$ with degree sequence $\pi$ has the Faber-Krahn property
  in the class $\mathcal{T}_\pi$ if and only if it is a \SLOst-tree.
  $G$ is then uniquely determined up to isomorphism.
\end{thm}

As an immediate corollary we get the result of \citet{Pruss:1998a}.

\begin{cor}[{\citep[Thm.~6.2]{Pruss:1998a}}]
  \label{thm:Tregular}
  In the class of semi-$d$-regular trees a graph $G$ has the
  Faber-Krahn property if and only if it is a \SLOst-tree.
  $G$ is then uniquely determined up to isomorphism.
\end{cor}

Before we proof these theorems we first want to show that each of
these two classes indead contains a \SLOst-tree.

\begin{lemma}
  \label{lem:unique-SLOst-tree}
  Each class $\mathcal{T}_\pi$ contains a \SLOst-tree that is uniquely
  determined up to isomorphism.
\end{lemma}
\begin{proof} 
  Let $\pi=(d_0,d_1,\ldots,d_{k-1},d_k,\ldots,d_{n-1})$ 
  be the degree sequence of $\mathcal{T}_\pi$, where 
  $2\leq d_0\leq d_1\leq\ldots\leq d_{k-1}$ and 
  $d_k = \ldots = d_{n-1} = 1$ (i.e., correspond to boundary
  vertices). 
  First we prove the existence of a \SLOst-tree by induction on
  $|\pi|$ (the number of vertices of $\pi$). 
  This is trivial for $|\pi|\leq 3$, since then $\pi=\{2,1,1\}$ and
  the corresponding graph is a path of length 2 and the vertex with
  degree 2 is choosen as root for the \SLO-ordering.

  Now we assume by induction that each $\mathcal{T}_\pi$ with 
  $|\pi|\leq n-1$ has a \SLOst-tree.
  For $|\pi|=n$ we construct a new degree sequence $\pi'$ by deleting
  the last $d_{k-1}-1$ elements from $\pi$ (which are all equal to $1$ as  
  $d_{k-1}-1<|\partial V| = n-k$) and by replacing $d_{k-1}$ by
  $d'_{k-1}=1$. Obviously $\pi'$ has $n-(d_{k-1}-1)$ elements. 
  By Lemma~\ref{lem:treesequence}, $\pi'$ is a tree sequence.
  By induction $\mathcal{T}_{\pi'}$ has a \SLOst-tree $T'$.
  Let $v$ be the first vertex of $T'$ w.\,r.\,t.\ the \SLO-ordering
  that is adjacent to some boundary vertex $w$.
  We replace $w$ by an interior vertex $u$ and add $d_{k-1}-1$ boundary
  vertices and get a tree $T$. Obviously $u$ has degree $d_{k-1}$ and thus
  $T$ has degree sequence $\pi$.
  Moreover, $T$ has a \SLOst-ordering which can be derived from the
  ordering in $T'$ by inserting the new vertex $u$ as the last
  interior vertex and the new boundary vertices as the last $d_{k-1}-1$
  vertices in the ordering. It is then easy to see that the properties
  (S1)--(S4) are satisfied.

  To show that two \SLOst-trees $G$ and $G'$ in a class
  $\mathcal{T}_\pi$ are isomorph we use a function $\phi$ that maps
  the vertex $v_i$ in the $i$-th position in the \SLOst-ordering of
  $G$ to the vertex $w_i$ in the $i$-th position in the
  \SLOst-ordering of $G'$.
  By the properities of the \SLOst-ordering, $\phi$ is an isomorphism,
  as $v_i$ and $w_i$ have the same degree and the images of all
  children of $v_i$ are exactly the children of $w_i$. The latter can
  be seen by looking on all interior vertices of $G$ in the reverse 
  \SLOst-ordering.
  Thus the proposition follows.
\end{proof} 
  

\section{Proof of the Theorems}

We first recall some basic results.
By definition the Laplace operator $\Delta$ is symmetric. Its
associate Rayleigh quotient on real valued functions $f$ on $V$ is the
fraction 
\begin{equation}
  \mathcal{R}_G(f)
  = \frac{\langle\Delta f,f\rangle}{\langle f,f \rangle}
  = \frac{\sum_{(u,v)\in E} (f(u)-f(v))^2}{\sum_{v\in V} f(v)^2}.
\end{equation}
For the Dirichlet operator $\Delta_0$ we get a similar Rayleigh
quotient. However, it is much simpler to consider $\mathcal{R}_G(f)$
again but restrict the set of functions $f$ such that $f(v)=0$ for all
boundary vertices $v\in\partial V$.
We denote the first Dirichlet eigenvalue of $\Delta_0(G)$ by
$\lambda(G)$. 
The following proposition states a well-known fact about Rayleigh
quotients.

\begin{proposition}
  \label{prop:mini-max-theorem}
  For a graph with boundary $G(V_0\cup\partial V, E_0\cup\partial E)$ 
  we have
  \begin{equation}
    \lambda(G) 
    = \min_{f\in\mathcal{S}} \mathcal{R}_G(f)
    = \min_{f\in\mathcal{S}} \frac{\langle\Delta f,f\rangle}{\langle f,f \rangle}
  \end{equation}
  where $\mathcal{S}$ is the set of all real-valued functions on $V$ with the
  constraint $f\big\vert_{\partial V}=0$.
  Moreover, if $\mathcal{R}_G(f)=\lambda(G)$ for a function $f\in\mathcal{S}$,
  then $f$ is an eigenfunction to the first Direchlet eigenvalue of
  $\Delta_0$.
\end{proposition}

For eigenfunctions of the Dirichlet operator the following
remarkable proposition holds.

\begin{proposition}
  Let $G(V_0\cup\partial V,E_0\cup\partial E)$ be a connected graph with
  boundary and $f$ an eigenfunction to some eigenvalue $\lambda$ of
  the Dirichlet operator. 
  Let $b_v$ denote the boundary vertices adjacent to $v$,
  i.e.\ $b_v=|\{w\in\partial V\colon (v,w)\in E\}|$.
  Then
  \begin{equation*}
    \lambda = \frac{\sum_{v\in V} b_v\,f(v)}{\sum_{v\in V} f(v)}\,.
  \end{equation*}
\end{proposition}
\begin{proof}
  Let $\mathbf{1}=(1,\ldots,1)'$ and 
  $i_v=|\{w\in V_0\colon (v,w)\in E\}|$ be the number of interior
  vertices adjacent to $v$. Thus $b_v+i_v=d_v$.
  A straight-forward computation gives
  \[
  \begin{split}
  \langle\mathbf{1},\Delta_0 f\rangle
  & \textstyle
    = \sum_{v\in V_0} d_v\,f(v) 
    - \sum_{v\in V_0}\sum_{\substack{(v,w)\in E\\ w\in V_0}} f(w) \\
  & \textstyle
    = \sum_{v\in V_0} d_v\,f(v) 
    - \sum_{w\in V_0} f(w)\sum_{\substack{(w,v)\in E\\ v\in V_0}} 1 \\
  & \textstyle
    = \sum_{v\in V_0} d_v\,f(v) - \sum_{w\in V_0} i_w\,f(w)
    = \sum_{v\in V_0} b_v\,f(v)\,.
  \end{split}
  \]
  Since $f$ is an eigenfunction we find 
  $\langle\mathbf{1},\Delta_0 f\rangle
  = \lambda\,\sum_{v\in V_0} f(v)$.
  As $f(v)=0$ for all boundary vertices $v\in\partial V$ the 
  result follows.
\end{proof}

\begin{proposition}[\citet{Friedman:1993a}]
  \label{prop:basic-properties}
  Let $G(V_0\cup\partial V,E_0\cup\partial E)$ be a connected graph with
  boundary.
  \begin{enumerate}[(1)]
  \item $\Delta_0(G)$ is a positive operator, i.e.\ $\lambda(G)>0$.
  \item An eigenfunction $f$ to the eigenvalue $\lambda(G)$
    is either positive or negative on all interior vertices of $G$.
  \item $\lambda(G)$ is monotone in $G$, i.e.\ if $G\subset G^\prime$ then
    $\lambda(G)>\lambda(G^\prime)$.
  \item $\lambda(G)$ is a simple eigenvalue.
  \end{enumerate}
\end{proposition}

\begin{remark}
  Let $T$ be a spannig tree of a graph $G$. By
  Prop.~\ref{prop:basic-properties} the first Dirichlet eigenvalue of
  the tree class $T\in\mathcal{T}$ is a lower bound for $\lambda(G)$.
\end{remark}

The main techniques for proving our theorems is \emph{rearranging} of
edges. We need two different types of rearrangement steps that we call
\emph{switching} and \emph{shifting}, respectively, in the following.

\begin{lemma}[{Switching, see also \citep[Lemma~5]{Leydold:1997a}}]
  \label{lem:switching}
  Let $G(V,E)$ be a tree with boundary in some
  $\mathcal{T}_\pi$.
  Let $(v_1,u_1),(v_2,u_2)\in E$ be edges such that $u_2$ is in the
  geodesic path from $v_1$ to $v_2$, but $u_1$ is not, see
  Fig.~\ref{fig:switching}.
  Then by replacing edges $(v_1,u_1)$ and $(v_2,u_2)$ by the 
  edges $(v_1,v_2)$ and $(u_1,u_2)$ we get a new 
  tree $G'(V,E')$ which is also contained in $\mathcal{T}_\pi$
  with the same set of boundary vertices.
  Moreover, we find for a function $f\in\mathcal{S}$
  \begin{equation}
    \label{eq:switching-inequality}
    \mathcal{R}_{G'}(f)\leq \mathcal{R}_{G}(f) 
  \end{equation}
  whenever $f(v_1)\geq f(u_2)$ and $f(v_2)\geq f(u_1)$.
  Inequality~(\ref{eq:switching-inequality}) is strict if both
  inequalities are strict.
\end{lemma}
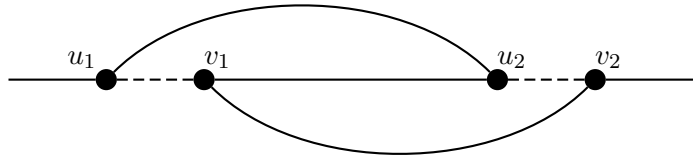
\begin{figure}[ht]
  \centering
  {
\psset{xunit=13mm}
\begin{pspicture}(0,-1)(7,1)
\pnode(0,0){L}
\vertex{U1}{1,0}\rput[rb](0.9,0.15){\small $u_1$}
\vertex{V1}{2,0}\rput[lb](2,0.15){\small $v_1$}
\vertex{U2}{5,0}\rput[lb](5,0.15){\small $u_2$}
\vertex{V2}{6,0}\rput[lb](6,0.15){\small $v_2$}
\pnode(7,0){R}
\edge{L}{U1}
\edgedash{U1}{V1}
\edge{V1}{U2}
\edgedash{U2}{V2}
\edge{V2}{R}
\nccurve[nodesep=0pt,angleA=45,angleB=135]{-}{U1}{U2}
\nccurve[nodesep=0pt,angleA=315,angleB=225]{-}{V1}{V2}
\end{pspicture}
}
  \caption{Switching: edges $(v_1,u_1)$ and $(v_2,u_2)$ are replaced
  by edges $(v_1,v_2)$ and $(u_1,u_2)$}
  \label{fig:switching}
\end{figure}
\begin{proof}
  Since by assumption $u_2$ is in the geodesic path from $v_1$ to
  $v_2$ and $u_1$ is not, $G'(V,E')$ again is a tree. The set of
  vertices does not change by construction. 
  Moreover, since this switching does not change the degrees of the
  vertices, the degree sequence remains unchanged.
  To verify inequality (\ref{eq:switching-inequality}) 
  we have to compute the effects of removing and inserting edges
  and get
  \[
  \begin{split}
    \langle\Delta(G') f,f\rangle - \langle\Delta(G) f,f\rangle
    &=
    \left[(f(v_1)-f(v_2))^2 + (f(u_1)-f(u_2))^2 \right] \\
    &\quad -\left[(f(v_1)-f(u_1))^2 + (f(v_2)-f(u_2))^2 \right] \\
    &=
    2\,(f(u_1)-f(v_2))\cdot (f(v_1)-f(u_2)) \\
    &\leq 0\,,
  \end{split}
  \]
  where last inequality is strict if both inequalities $f(v_1)\geq
  f(u_2)$ and $f(v_2)\geq f(u_1)$ are strict.
  Thus the proposition follows.
\end{proof}

\begin{cor}
  \label{lem:switching-ef}
  Let $G(V,E)$ be a tree with boundary in some
  $\mathcal{T}_\pi$ and let $G'(V,E')$ be a tree obtained from $G$
  by applying Switching as defined in Lemma~\ref{lem:switching}.
  If $f$ is a non-negative eigenfunction to the first Dirichlet
  eigenvalue of $G$ then $\lambda(G')\leq\lambda(G)$
  whenever $f(v_1)\geq f(u_2)$ and $f(v_2)\geq f(u_1)$.
  Moreover, $\lambda(G')<\lambda(G)$ if one of these two
  inequalities is strict.
\end{cor}
\begin{proof}
  The first inequality is an immediate consequence of
  Lemma~\ref{lem:switching} and Prop.~\ref{prop:mini-max-theorem}
  \[
  \lambda(G')\leq\mathcal{R}_{G'}(f)\leq\mathcal{R}_{G}(f) = \lambda(G).
  \]
  For the second statement notice that $\lambda(G')=\lambda(G)$ if and
  only if $\mathcal{R}_{G'}(f)=\mathcal{R}_{G}(f)$ and $f$ is an
  eigenfunction to $\lambda(G')$ on $G'$, since $\lambda(G')$ is simple
  (Props.~\ref{prop:mini-max-theorem} and \ref{prop:basic-properties}).
  Therefore, if $\lambda(G')=\lambda(G)$ we find
  \[
  \begin{split}
    \lambda(G)f(v_1) &= \Delta(G) f(v_1)
    = d_{v_1}f(v_1) - f(u_1) - \sum_{\substack{(v_1,w)\in E\\ w\not=u_1}} f(w) \\
    = \lambda(G')f(v_1) &= \Delta(G') f(v_1)
    = d_{v_1}f(v_1) - f(v_2) - \sum_{\substack{(v_1,w)\in E'\\ w\not=v_2}} f(w)\,. \\
  \end{split}
  \]
  Since the summation is done over the same neighbors of $v_1$ in this
  equation we find $f(u_1) = f(v_2)$. Analogously we derive from
  $\Delta(G)f(u_1)=\Delta(G')f(u_1)$,
  $f(v_1)=f(u_2)$. Thus the proposition follows.
\end{proof}
\begin{remark}
  Lemma~\ref{lem:switching} and Cor.~\ref{lem:switching-ef} hold
  analogously for arbitrary graphs.
\end{remark}

\begin{lemma}[Shifting]
  \label{lem:shifting}
  Let $G(V,E)$ be a tree with boundary in some graph class
  $\mathcal{T}$.
  Let $(u,v_1)\in E$ be an edge and $v_2\in V$ some vertex such that
  $u$ is not in the geodesic path from $v_1$ to $v_2$,
  Fig.~\ref{fig:shifting}. 
  Then by replacing edge $(u,v_1)$ by the edge $(u,v_2)$ we get a new 
  tree $G'(V,E')$ which is also contained in $\mathcal{T}$.
  If $v_2\in V_0$ is an interior vertex then the number of boundary
  vertices remains unchanged.
  Moreover, we find for a non-negative function $f\in\mathcal{S}$
  \begin{equation}
    \label{eq:shifting-inequality}
    \mathcal{R}_{G'}(f)\leq \mathcal{R}_{G}(f) 
  \end{equation}
  if and only if $f(v_1)\geq f(v_2)$.
  The inequality is strict if $f(v_1) > f(v_2)$.
\end{lemma}
Notice that if $G$ is in some class $\mathcal{T}_d$ (or
$\mathcal{T}_\pi$) then in general $G'$ need not be a member of
this graph class any more.
\begin{figure}[ht]
  \centering
  {
\psset{xunit=13mm}
\begin{pspicture}(0,-2)(4,0.5)
\pnode(0,0){L}
\vertex{V1}{1,0}\rput[lb](1.05,0.15){\small $v_1$}
\vertex{V2}{3,0}\rput[lb](3.05,0.15){\small $v_2$}
\vertex{U}{1,-1}\rput[rt](0.9,-1.15){\small $u$}
\pnode(4,0){R}
\pnode(1,-2){D}
\edge{L}{V1}
\edge{V1}{V2}
\edgedash{U}{V1}
\edge{U}{V2}
\edge{V2}{R}
\edge{U}{D}
\end{pspicture}
}
  \caption{Shifting: edge $(u,v_1)$ is replaced by edge $(u,v_2)$}
  \label{fig:shifting}
\end{figure}
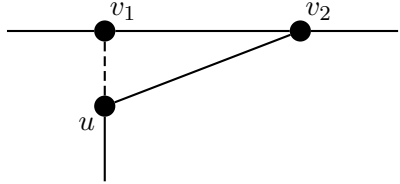
\begin{proof}
  Analogously to the proof of Lemma~\ref{lem:switching}.
\end{proof}
\begin{remark}
  Lemma~\ref{lem:shifting} holds analogously for arbitrary graphs.
\end{remark}

We now can use a sequence of switchings and shiftings to transform
any tree $G$ with boundary in some class $\mathcal{T}_\pi$ into
\SLOst-tree $G^\ast\in\mathcal{T}_\pi$.

\begin{lemma}
  \label{lem:SLO-tree}%
  \label{lem:necessary-SLO-tree}%
  \label{lem:ef-decreasing}%
  Let $G(V,E)$ be a tree with boundary in some class
  $\mathcal{T}_\pi$.
  Then there exists a \SLO-tree $G'(V,E')$
  in $\mathcal{T}_\pi$ with $\lambda(G')\leq\lambda(G)$.

  Furthermore, if $G$ has the Faber-Krahn property then there exists
  already a \SLO-ordering $\prec$ of the vertices (i.e., $G$ is a
  \SLO-tree). If, moreover, $f$ is a non-negative eigenfunction to
  $\lambda(G)$ then $v\prec w$ implies $f(v)\geq f(w)$.
\end{lemma}
\begin{proof}
  Let $n=|V|$ and $k=|V_0|$ denote the number of vertices and of
  interior vertices of $G$, respectively, and let $f$ be a
  non-negative eigenfunction to the first Dirichlet eigenvalue of $G$.
  We assume that the vertices of $G$, 
  $V=\{v_0,v_1,\ldots,v_{k-1},v_k,\ldots,v_{n-1}\}$,
  are numbered such that $f(v_i)\geq f(v_j)$ if $i<j$, i.e., they are
  sorted with respect to $f(v)$ in non-increasing order. 
  We define a well-ordering $\prec$ on $V$ by
  $v_i\prec v_j$ if and only if $i<j$.

  Now we use a series of switchings to construct the desired
  new tree $G'$. This is done recursively such that we have a ball 
  that has already the desired \SLO-ordering in the center of each
  intermediate graph. This ball grows in every recursion step until
  all vertices of the initial graph $G$ are used.

  We start with the first vertex $v_0$ of this ordered set of
  vertices. If $v_0$ is adjacent to $v_1$ there is nothing to do.
  Else, we check whether $v_0$ is adjacent to some vertex 
  $w$ with $f(w)=f(v_1)$ and $v_1\prec w$. If there exists such a
  vertex we just exchange the positions of these two vertices in the
  ordering of $V$ (and update the indices of the vertices).
  (In particular this is the case when $v_1$ is a boundary vertex then 
  by our assumptions $0\leq f(w)\leq f(v_1)=0$ and thus
  $f(w)=f(v_1)=0$ and this condition is satisfied.)
  Otherwise, there exists a child vertex $u_0$ of $v_0$
  with $v_1\prec u_0$ and a path $P_{0,1}$ from $v_0$ to $v_1$, since
  $G$ is connected.
  There also exist a parent of $v_1$ (which is in this path
  $P_{0,1}$ and which cannot be $v_0$) and some child vertices (which
  are not in this path).
  The latter exist as $v_1$ cannot be a boundary vertex, since one of
  the above two cases would apply.
  Now if $u_0 \in P_{0,1}$ then let $u_1$ be one these child
  vertices; else let $u_1$ we the parent of $v_1$.
  As by construction $v_0\prec v_1 \prec u_0, u_1$ we have $f(v_0)\geq
  f(v_1) \geq f(u_0), f(u_1)$ and hence we can apply
  Lemma~\ref{lem:switching}, exchange edges $(v_0,u_0)$ and
  $(v_1,u_1)$ by $(v_0,v_1)$ and $(u_0,u_1)$, 
  and get a new graph $G_1$ with 
  $\mathcal{R}_{G_1}(f)\leq\mathcal{R}_{G}(f)$
  which also belongs to $\mathcal{T}_\pi$.

  By this switching step we have exchanged a child of $v_0$
  by $v_1$ (if necessary) which then becomes a child of $v_0$.
  By the same procedure we can exchange all
  other vertices adjacent to $v_0$ with the respective vertices $v_2$,
  $v_3$, \ldots, $v_{s_0}$, where $s_0=d_{v_0}$, and get graphs $G_2$,
  $G_3$, \ldots, $G_{s_0}$ in $\mathcal{T}_\pi$ with  
  $\mathcal{R}_{G_i}(f)\leq \mathcal{R}_{G_{i-1}}(f)$.

  Next we proceed in an analogous manner with all children $u$
  of $v_1$ with $v_1\prec u$ and make all vertices $v_{s_0+1}$,
  $v_{s_0+2}$, \ldots, $v_{s_1}$ adjacent to $v_1$, where $s_1=s_0+d_{v_1}-1$,
  and get graphs $G_{s_0+1}$, $G_{s_0+2}$, \ldots, $G_{s_1}$.
  By processing all interior vertices in this way we get a sequence of
  graphs 
  \begin{equation}
    \label{eq:rearrange-sequence-SLO-tree}
    G = G_0 \rightarrow G_1 \rightarrow G_2 \rightarrow \ldots
    \rightarrow G_k = G'
  \end{equation}
  in $\mathcal{T}_\pi$ with
  \begin{equation}
    \label{eq:rearrange-inequ-SLO-tree}
    \lambda(G) = \mathcal{R}_{G_0}(f)\geq
    \mathcal{R}_{G_1}(f)\geq \ldots \geq
    \mathcal{R}_{G_k}(f)\geq \lambda(G')\;.
  \end{equation}
  In step $G_{r-1}\rightarrow G_r$ there is either nothing to do (when
  we assume that the vertices are already in the proper ordering),
  or the vertex $v_r$ is made adjacent to vertex $v_{r-1}\prec v_r$ by
  a switching step: 
  Let $P_{r-1,r}$ be the geodesic path from $v_{r-1}$ to $v_{r}$. 
  By construction of our sequence of graphs we have 
  $h(v_{r-1}) \leq h(v_r)$ in graph $G_{r-1}$ and thus the parent
  $w_r$ of $v_r$ must be in $P_{r-1,r}$.
  Moreover, $v_r$ cannot be a boundary vertex (since otherwise we can
  use the argument from above and we only had to change the ordering
  of the vertices) and thus has some child $u_r$.
  Furthermore this path either contains some child $u_{r-1}$ of
  $v_{r-1}$, or it contains the parent of $v_{r-1}$. In the latter
  case there exists at least one child $u_{r-1}$.
  Now we can use switching and replace either
  edges $(v_{r-1},u_{r-1})$ and $(v_r,u_r)$ by the edges 
  $(v_{r-1},v_r)$ and $(u_{r-1},u_r)$ (if $u_{r-1}$ is contained in
  $P_{r-1,r}$) 
  or (otherwise) edges $(v_{r-1},u_{r-1})$ and $(w_r,v_r)$ by the
  edges $(v_{r-1},v_r)$ and $(u_{r-1},w_r)$.
  In both cases we have can apply Lemma~\ref{lem:switching} as
  $f(v_{r-1})\geq f(v_r)\geq f(u_{r-1}),\,f(w_r)\,,f(u_r)$.
  (It cannot happen that $v_r$ is adjacent to some vertex $w$ with
  $w\prec v_{r-1}$.)
  In the consecutive steps edges between vertices $u$ and $w$ with
  $u\prec w\prec v_{r+1}$ are neither deleted nor inserted any more.
  Hence 
  $\lambda(G')\leq\mathcal{R}_{G'}(f)\leq \mathcal{R}_{G}(f)=\lambda(G)$.

  It remains to show that $\prec$ is a \SLO-ordering of the vertices
  $V$ in $G'$. Property (S3) holds by definition of the ordering
  $\prec$. By construction (S2) holds. Moreover, $G'$ is
  built by stepwise adding layers to a ball. Thus property (S1) holds
  and the first statement follows.

  Now assume that $G$ has the Faber-Krahn property. Then equality
  holds in (\ref{eq:rearrange-inequ-SLO-tree}) everywhere.
  Furthermore, $f$ must be an eigenfunction to the first Dirichlet
  eigenvalue for every graph $G_i$ in this sequence. Otherwise, if $f$
  is not an eigenfunction of a graph $G_i$ 
  then $\lambda(G_i)<\mathcal{R}_{G_i}(f) = \lambda(G)$, by
  Prop.~\ref{prop:mini-max-theorem} a contradiction.

  For switching step $G_{r-1}\rightarrow G_r$ we have
  $f(v_r)\geq f(u_{r-1})$. If $f(v_r)=f(u_{r-1})$ there would be
  nothing to do (we only change the positions of $v_r$ and $u_{r-1}$
  in the ordering $\prec$). Hence we have $f(v_r)>f(u_{r-1})$ and by 
  Corollary~\ref{lem:switching-ef},
  $\lambda(G_r)<\lambda(G_{r-1})$, a contradiction to the Faber-Krahn
  property of $G$.

  The monotinicity property of $f$ follows by the same reasons.
\end{proof}

\begin{lemma}
  \label{lem:SLOst-tree}
  Let $G(V,E)$ be a tree with boundary in some
  $\mathcal{T}_\pi$.
  Then there exists a \SLOst-tree $G^\ast(V,E^\ast)$
  in $\mathcal{T}_\pi$ with $\lambda(G^\ast)\leq\lambda(G)$.
\end{lemma}
\begin{proof}
  Let again $n=|V|$ and $k=|V_0|$ denote the number of vertices and of
  interior vertices of $G$, respectively, and let $f$ be a
  non-negative eigenfunction to the first Dirichlet eigenvalue of $G$.
  Then by Lemma~\ref{lem:SLO-tree} there exists a \SLO-tree
  $G_0'=G'(V,E')$ in $\mathcal{T}_\pi$ with the \SLO-ordering $\prec$.
  The vertices of $G$ (and $G'$) 
  $V=\{v_0,v_1,\ldots,v_{k-1},v_k,\ldots,v_{n-1}\}$ are numbered
  such that $v_i\prec v_j$ if and only if $i<j$.
  Moreover, by the construction in the proof of
  Lemma~\ref{lem:SLO-tree} we find $f(v)\geq f(w)$ if $v\prec w$.
  The degree sequence of $G$ is given by
  $\pi=(d_0,d_1,\ldots,d_{k-1},d_k,\ldots,d_{n-1})$ 
  such that the degrees $d_i$ are non-decreasing for $0\leq i<k$, and
  $d_j=1$ for $j\geq k$ (i.e., correspond to boundary vertices).

  Now we start with root $v_0$. If $d_{v_0}=d_0$
  ($=\min_{0\leq i\leq k} d_i$) then there is nothing to do.
  Otherwise, we can use shifting to replace all edges
  $(v_0,v_{d_0+1})$, $(v_0,v_{d_0+2})$, \ldots, by the respective 
  edges $(v_1,v_{d_0+1})$, $(v_1,v_{d_0+2})$, \ldots\,.
  As $v_0\prec v_1$ we have $f(v_0)\geq f(v_1)$ and thus we can apply 
  Lemma~\ref{lem:shifting} and get a new graph $G'_1$ with 
  $\mathcal{R}_{G'_1}(f)\leq \mathcal{R}_{G'}(f)$.
  Notice that $G'_1$ is again a \SLO-tree. However, it might happen
  that the degree sequence has changed and
  $G'_1\not\in\mathcal{T}_\pi$.

  Next we proceed in the same way with vertex $v_1$.
  We denote the degree of a vertex $v_j$ in a graph $G_i'$ with index
  $i$ by $d^{(i)}_{v_j}$.
  Notice that 
  $d^{(1)}_{v_1}\geq \min_{1\leq i\leq k} d_i=d_1$.
  If $d^{(1)}_{v_1}=d_1$ there is nothing to do.
  Otherwise, we can use shifting to replace all edges 
  $(v_1,v_{s_1+1})$, $(v_1,v_{s_1+2})$, \ldots, by the respective
  edges $(v_2,v_{s_1+1})$, $(v_2,v_{s_1+2})$, \ldots,
  where $s_1=d_0+d_1$. Again we can apply Lemma~\ref{lem:shifting}
  and get a new graph $G'_2$ with 
  $\mathcal{R}_{G'_2}(f)\leq \mathcal{R}_{G'_1}(f)$.
  We can continue in this way and get a sequence of \SLO-trees
  \begin{equation}
    \label{eq:rearrange-sequence-SLOst-tree}
    G \rightarrow G' = G'_0 \rightarrow G'_1 \rightarrow G'_2
    \rightarrow \ldots \rightarrow G'_k = G^\ast
  \end{equation}
  with
  \begin{equation}
    \label{eq:rearrange-inequ-SLOst-tree}
    \lambda(G) = \mathcal{R}_{G}(f)\geq 
    \mathcal{R}_{G'_0}(f)\geq 
    \mathcal{R}_{G'_1}(f)\geq \ldots \geq
    \mathcal{R}_{G'_k}(f)\geq \lambda(G^\ast)\;.
  \end{equation}
  Notice that we always have $d^{(r)}_{v_r}\geq d_r$. This follows
  from the fact that 
  $\sum_{j\leq r} d^{(0)}_{v_j} \geq \sum_{j\leq r} d_j$
  as the right hand side of this inequality is the minimum of any sum
  of degrees of $j$ interior vertices of $G'$.
  Morever by our construction,
  $\sum_{j\leq r} d^{(r)}_{v_j} = \sum_{j\leq r} d^{(0)}_{v_j}$
  and $\sum_{j<r} d^{(r)}_{v_j} = \sum_{j<r} d_j$.
  Hence 
  $d^{(r)}_{v_r} 
  = \sum_{j\leq r} d^{(r)}_{v_j} - \sum_{j<r} d^{(r)}_{v_j}
  = \sum_{j\leq r} d^{(0)}_{v_j} - \sum_{j<r} d_j
  \geq \sum_{j\leq r} d_j - \sum_{j<r} d_j
  = d_r$.
  In step $G'_r\rightarrow G'_{r+1}$ there is either nothing to do,
  or edges are exchanged such that vertex $v_r$ has the desired
  degree. In the consecutive steps edges that are incident to a vertex
  $u\prec v_{r+1}$ are neither deleted nor inserted. 

  The resulting \SLO-tree $G^\ast$ has the same degree sequence $\pi$
  as $G$ and thus belongs to class $\mathcal{T}_\pi$. It also
  satisfies property $(S4)$, i.e.\ $\prec$ is a \SLOst-ordering of the
  vertices.
\end{proof}

For our theorem on the class $\mathcal{T}_d$ we need a modified
version of this lemma. To state this new proposition we need a partial
ordering of degree sequences.
Let $\pi=(d_0,d_1,\ldots,d_{k-1},d_k,\ldots,d_{n-1})$ and 
$\pi'=(d'_0,d'_1,\ldots,d'_{k'-1},d_{k'},\ldots,d_{n-1})$ be two
degree sequence of some trees with the same number of vertices $n$
and respective numbers $k$ and $k'$ of interior vertices (not
necessarily equal). Again we assume that the first $k$ (and $k'$,
resp.)\ degrees correspond to the interior vertices and are ordered
non-decreasingly. 
Then we write $\pi\unlhd\pi'$ if the above condition holds and
$\sum_{j\leq r} d_j\leq\sum_{j\leq r} d'_j$ for all $0\leq r<n$.
\begin{lemma}
  \label{lem:SLOst-tree-d}
  Let $G(V,E)$ be a tree with boundary with degree sequence $\pi$ and
  let $\pi'$ another degree sequence with $\pi'\unlhd\pi$.
  Then there exists a \SLOst-tree $G^\ast(V,E^\ast)$
  in $\mathcal{T}_{\pi'}$ with $\lambda(G^\ast)\leq\lambda(G)$.
\end{lemma}
\begin{proof}
  Completely analogous to the proof of Lemma~\ref{lem:SLOst-tree}.
\end{proof}
Notice that Lemma~\ref{lem:SLOst-tree} is a special case of this lemma
as $\pi\unlhd\pi$.
It can also be applied to prove Theorem~\ref{thm:Tdeg} for class
$\mathcal{T}_d$ as we immediately have $\pi^\circ\unlhd\pi$ with
$\pi^\circ=(d,d,\ldots,d,d^\circ,1,\ldots,1)$
where $d^\circ=d+\sum_{v\in V_0}(d_v-d)$.

Next we show that every tree with the Faber-Krahn property has a
\SLOst-ordering.

\begin{lemma}
  \label{lem:decreasing-child}
  Let $G$ be a \SLO-tree with a non-negative eigenfunction $f$ of
  $\lambda(G)$. Then every interior vertex $v$ has a child $w$ with
  $f(w)<f(v)$. 
\end{lemma}
\begin{proof}
  First assume $v$ that is not the root of $G$. 
  Let $u$ be the parent of $v$. Then by
  Lemma~\ref{lem:ef-decreasing} $f(v)\leq f(u)$ and
  $f(v)\geq f(w)$ for all children $w$ of $v$.
  Now suppose that $f(v)=f(w)$ for all children of $v$.
  Then
  $\lambda(G) f(v) 
  = \Delta f(v)
  = \sum_{(v,x)\in E} (f(v)-f(x))
  = f(v)-f(u)\leq 0$,
  a contradiction as both $f(v)>0$ and $\lambda(G)>0$ by
  Prop.~\ref{prop:basic-properties}. 
  If $v$ is the root of $G$ then all vertices adjacent to $v$ are
  children of $v$. If we again suppose for all these children we
  $f(w)=f(v)$ then we find analogously 
  $\lambda(G) f(v) = 0$, again a contradiction.
\end{proof}

\begin{lemma}
  \label{lem:isomorphic-subtrees} 
  Let $G(V,E)$ be a \SLOst-tree and $f$ a non-negative eigenfunction
  to $\lambda(G)$. Let $v$ and $w$ two vertices $v$ and $w$ with 
  $f(v)=f(w)$. Then the subtrees $T_v$ and 
  $T_w$ rooted at $v$ and $w$, respectively, are isomorph.
\end{lemma}
\begin{proof}
  We prove this lemma by induction from boundary vertices to the
  root $v_0$. It is obviously trivial for boundary vertices.
  Without loss of generality we assume $v\prec w$.

  We start with the case where $v$ is not the root $v_0$ of
  \SLOst-ordering.
  Let $u_v$ and  $u_w$ be the parents of $v$ and $w$,
  respectively. Then from $\Delta(G)f(v)$ and $\Delta(G)f(w)$ we get
  $f(u_v)=(d_v-\lambda(G))\,f(v)
  -\sum_{\substack{(v,x)\in E\\ x\not=u_v}} f(x)$
  and
  $f(u_w)=(d_w-\lambda(G))\,f(w)
  -\sum_{\substack{(w,y)\in E\\ y\not=u_w}} f(y)$.
  By property (S2) and Lemma~\ref{lem:ef-decreasing}
  we have $f(u_v)\geq f(u_w)$ and therefore it follows from
  $f(v)=f(w)$,
  \begin{equation}
    \label{eq:isomorphic-subtrees-1}
    (d_{w}-d_{v})\,f(v) \leq
    \sum_{\substack{(w,y)\in E\\ y\not=u_w}} f(y)
    -\sum_{\substack{(v,x)\in E\\ x\not=u_v}} f(x)
  \end{equation}
  where the sums on the right hand side are over all children of $w$
  and $v$, respectively. Let $m$ be a child of $v$ such that 
  $f(m)\leq f(x)$ for all children $x$ of $v$. Notice that by (S2)
  $x\prec y$ and thus by Lemma~\ref{lem:ef-decreasing}
  $f(x)\geq f(y)$ for all children $y$ of $w$;
  in particular $f(m)\geq f(y)$. Thus
  $\sum_{\substack{(v,x)\in E\\ x\not=u_v}} f(x)\geq (d_v-1)\,f(m)$
  and
  $\sum_{\substack{(w,y)\in E\\ y\not=u_w}} f(y)\leq (d_w-1)\,f(m)$.
  Consequently
  \begin{equation}
    \label{eq:isomorphic-subtrees-2}
    \sum_{\substack{(w,y)\in E\\ y\not=u_w}} f(y)
    -\sum_{\substack{(v,x)\in E\\ x\not=u_v}} f(x)
    \leq (d_{w}-d_{v})\,f(m)
  \end{equation}
  and by (\ref{eq:isomorphic-subtrees-1}) 
  $(d_{w}-d_{v})\,f(v) \leq (d_{w}-d_{v})\,f(m)$.

  By Prop.~\ref{prop:basic-properties} and
  Lemma~\ref{lem:decreasing-child}, $0<f(m)<f(v)$.
  By property (S4), $d_{v}\leq d_{w}$.
  Hence $d_{v}=d_{w}$. Then the right hand side of
  (\ref{eq:isomorphic-subtrees-1}) 
  (and left hand side of (\ref{eq:isomorphic-subtrees-2}))
  vanishes and $f$ must have the same value for all children of $v$
  and $w$ (in particular $f(x)=f(y)$).
  It then follows by induction that $T_v$ and $T_w$ are isomorph.

  The case where $v$ is the root $v_0$ of \SLOst-ordering,
  remains. Then we set $u_v = v_1$ and all estimations are still
  valid. Thus the proposition follows.
\end{proof}

\begin{lemma}
  \label{lem:necessary-SLOst-tree}
  If a tree $G(V,E)$ with boundary has the Faber-Krahn
  property in some class $\mathcal{T}_\pi$, then $G$ is a \SLOst-tree.
\end{lemma}
\begin{proof}
  By Lemma~\ref{lem:necessary-SLO-tree} $G$ is a \SLO-tree.
  In the proof of Lemma~\ref{lem:SLOst-tree} we have produced the
  sequence (\ref{eq:rearrange-sequence-SLOst-tree}) of trees
  where the inequalities (\ref{eq:rearrange-inequ-SLOst-tree}) hold.
  Since $G$ has the Faber-Krahn property, equality holds in each
  of these inequalities.
  Notice that $G'$ and $G^\ast$ are in class $\mathcal{T}_\pi$ while
  all other graphs $G'_i$ need not. However, for every graph $G'_i$ in
  this sequence that belongs to $\mathcal{T}_\pi$ we have by the
  Faber-Krahn property $\lambda(G'_i)=\lambda(G)$ and $f$ is also an
  eigenfunction to the first Dirichlet eigenvalue of $G_i$.
  Otherwise we had $\lambda(G_i)<\mathcal{R}_{G_i}(f)=\lambda(G)$, 
  a contradiction.

  Now suppose there is a graph $G_r\in\mathcal{T}_\pi$ while 
  $G_{r+1}\not\in\mathcal{T}_\pi$.
  We denote the children of vertex $v_r$ in $G_r$ by $w_1,\ldots,w_s$
  and its parent by $u_r$. In step $G'_r\rightarrow G'_{r+1}$ we
  replace the edges $(v_r,w_{d_r})$, \ldots, $(v_r,w_s)$ by the
  respective edges $(v_{r+1},w_{d_r})$, \ldots, $(v_{r+1},w_s)$.
  Hence $s>d_r-1$, since otherwise there would be nothing to do and
  $G_{r+1}=G_r$, a contradiction to $G_{r+1}\not\in\mathcal{T}_\pi$.
  Notice that the neighbors of $v_r$ in $G_{r+1}$ do not change any
  more in the subsequent steps.
  As $f$ is an eigenfunction to both $G_r$ and $G^\ast$ to the same
  eigenvalue $\lambda(G)$ it follows that 
  $\Delta(G'_r) f(v_r) = \Delta(G^\ast) f(v_r)$, i.e.\
  \[
  (s+1) f(v_r) - f(u_r) -\sum_{j=1}^s f(w_j)
  = d_r f(v_r) - f(u_r) -\sum_{j=1}^{d_r-1} f(w_j)
  \]
  and thus $(s-d_r+1) f(v_r)=\sum_{j=d_r}^s f(w_j)$.
  Since $f(v_r)\geq f(w_1)\geq f(w_j)\geq f(w_s)\geq 0$
  for all $j=1,\ldots,s$ by Lemma~\ref{lem:ef-decreasing}, we find
  $f(v_r)=f(w_j)$ for all children $w_j$, a contradiction to
  Lemma~\ref{lem:decreasing-child}.
  If $r=0$, i.e. $v_r$ is the root and there is no parent of $v_r$,
  then same argment and holds analogously.

  Hence there cannot be a graph $G_r\in\mathcal{T}_\pi$ while 
  $G_{r+1}\not\in\mathcal{T}_\pi$. Therefore each graph $G'_i$ in
  sequence (\ref{eq:rearrange-sequence-SLOst-tree}) belongs to class
  $\mathcal{T}_\pi$ and $f$ is an eigenfunction for each of these.
  We show for each $r$ that $G_r$ is isomorph to $G_{r+1}$ and
  consequently isomorph to $G^\ast$. Thus all these graphs, in
  particular $G'_0$, are \SLOst-trees.
  Notice that for step $G'_r\rightarrow G'_{r+1}$ we either find
  $G_r=G_{r+1}$, or $f(v_r)=f(v_{r+1})$, since otherwise we had
  $\mathcal{R}_{G'_r}(f)>\mathcal{R}_{G'_{r+1}}(f)$ by
  Lemma~\ref{lem:shifting}.
  In the first case there remains nothing to show. In the latter case
  the subtrees (of both $G_r$ and $G_{r+1}$) rooted at the
  respective vertices $v_r$ and $v_{r+1}$ are isomorphic by
  Lemma~\ref{lem:isomorphic-subtrees}.  
  As only edges incident to $v_r$ are shifted to $v_{r+1}$ the
  isomorphism between $G_r$ and $G_{r+1}$ follows.
\end{proof}

Now we are ready to prove our theorems.

\begin{proof}[Proof of Theorem~\ref{thm:Tdegseq}]
  The necessity of the condition has been shown in
  Lemma~\ref{lem:necessary-SLOst-tree}. The sufficiency follows from
  the fact that \SLOst-trees are unique determined up to isomorphism
  (Lemma~\ref{lem:unique-SLOst-tree}).
\end{proof}

\begin{proof}[Proof of Theorem~\ref{thm:Tdeg}]
  Let $\pi=(d_0,d_1,\ldots,d_{k-1},1,\ldots,1)$ be the degree
  sequence of $G$, where $d\leq d_0\leq d_1\leq\ldots\leq d_{k-1}$ are
  the degrees for the interior vertices. 
  Define a new degree sequence by
  $\pi^\circ=(d,d,\ldots,d,d^\circ,1,$ $\ldots,1)$
  where $d^\circ=d+\sum_{v\in V_0}(d_v-d)$.
  Then $\pi'\unlhd\pi$ and we can apply Lemma~\ref{lem:SLOst-tree-d}.
  The necessity of the condition follows analogously to the proof
  Lemma~\ref{lem:necessary-SLOst-tree}.
  The sufficiency follows from the fact that \SLOst-trees are unique
  determined up to isomorphism (Lemma~\ref{lem:unique-SLOst-tree}).
\end{proof}

\begin{proof}[Proof of Theorem~\ref{thm:klobuerstel}]
  This is an immediate corollary of Thm.~\ref{thm:Tdeg} as
  $\mathcal{T} = \mathcal{T}_2$.
\end{proof}


\section{Further Results}
\label{sec:further-results}

One might ask what happens when we relax the conditions in the class
$\mathcal{T}^{(n,k)}$ and $\mathcal{T}_d^{(n,k)}$.
We then get the following classes
\begin{align}
  \mathcal{T}^{(n,\cdot)} 
  &= \{\text{$G$ is a tree, with $|V|=n$}\} \\
  \mathcal{T}_d^{(n,\cdot)} 
  &= \{G\in\mathcal{T}^{(n,\cdot)}\colon d_v\geq d \text{ for all }v\in V_0\}
\end{align}
where we keep the total number of vertices fixed, and
\begin{align}
  \mathcal{T}^{(\cdot,k)} 
  &= \{\text{$G$ is a tree, with $|V_0|=k$}\} \\
  \mathcal{T}_d^{(\cdot,k)} 
  &= \{G\in\mathcal{T}^{(\cdot,k)}\colon d_v\geq d \text{ for all }v\in V_0\}
\end{align}
where we keep the number of interior vertices fixed.
Using the arguments from the proofs of our theorems we find the
following characterizations for graphs with the Faber-Krahn property.

\begin{thm}
  A tree $G$ with boundary has the Faber-Krahn property
  \begin{enumerate}[(i)]
  \item 
    in $\mathcal{T}^{(n,\cdot)}$ if and only if it is a 
    path with $n$ vertices.
    (This is the result of \citep{Katsuda;Urakawa:1999a}.)
  \item 
    in $\mathcal{T}_d^{(n,\cdot)}$ if and only if it is a \SLOst-tree
    where exactly one interior vertex has degree 
    $d^\circ$ with $d\leq d^\circ<2\,d$
    and all other interior vertices have degree $d$.
    (This is the \SLOst-tree in $\mathcal{T}_d^{(n,\cdot)}$
    with the greatest number of interior vertices.)
  \item 
    in $\mathcal{T}^{(\cdot,k)}$ if and only if it is a 
    path with $k+2$ vertices.
  \item 
    in $\mathcal{T}_d^{(\cdot,k)}$ if and only if it is a \SLOst-tree
    where all interior vertices have degree $d$.
  \end{enumerate}
  $G$ is then uniquely determined up to isomorphism.
\end{thm}

For the classes $\mathcal{T}_\pi$ we cannot give a similar theorem. 
However, we can ask whether we can compare the least first Dirichlet
eigenvalue in classes with the same number of vertices. From
Lemma~\ref{lem:SLOst-tree-d} we can derive the following result.

\begin{thm}
  Let $\pi$ and $\pi'$ be two tree sequences with $|\pi|=|\pi'|$ and
  let $G$ and $G'$ be trees with the Faber-Krahn property in
  $\mathcal{T}_\pi$ and $\mathcal{T}_{\pi'}$, respectively.
  If $\pi'\unlhd\pi$ then $\lambda(G)\leq\lambda(G')$ where equality
  holds if and only if $\pi=\pi'$.
\end{thm}


\section*{Acknowledgement}

The authors would like to thank Franziska Berger for helpful
discussions.


\bibliographystyle{abbrvnat}
\bibliography{local}


\end{document}